\documentclass[12pt]{amsart}
\usepackage{amssymb,amsmath}

\begin{document}

\flushbottom

\title{Special It\^o maps and an $L^2$ 
Hodge theory for one forms on path spaces}
\author{K. D. Elworthy }
\address{K. D. Elworthy\\
MATHEMATICS INSTITUTE, UNIVERSITY OF WARWICK, COVENTRY CV4 7AL, UK}
\author{Xue-Mei Li}
\address{ Xue-Mei Li\\
DEPARTMENT OF MATHEMATICS, UNIVERSITY OF CONNECTICUT, 196 AUDITORIUM ROAD,
STORRS, CT06269-3009, USA.}
 \email{xmli@math.uconn.edu}   

\date{}     \maketitle

\newcommand{\A}{{\bf \mathcal A}}
\newcommand{\B}{{ \bf \mathcal B }}
\newcommand{\C}{{\mathcal C}}
\newcommand{\D}{{\rm I \! D}}
\newcommand{\E}{{\Bbb E}}
\newcommand{\F}{{\mathcal F}}
\newcommand{\G}{{\mathcal G}}
\newcommand{\K}{{\mathcal K}}
\newcommand{\p}{{\Bbb P}}
\newcommand{\R}{{\Bbb R}}
\newcommand{\I}{{\mathcal I}}

\newtheorem{theorem}{Theorem}[section]
\newtheorem{proposition}[theorem]{Proposition} 
\newtheorem{lemma}[theorem]{Lemma} 
\newtheorem{corollary}[theorem]{Corollary} 
\newtheorem{definition}{Definition}[section]

\renewcommand{\thefootnote}{}

\footnote{Research partially supported by NSF, EPSRC GR/NOO 845,
  the Alexander von Humboldt stiftung, and
EU grant ERB-FMRX-CT96-0075.}
\footnote{This article is published in `Stochastic processes, physics and geometry: new interplays, I' (Leipzig, 1999), 145–162, CMS Conf. Proc., 28, Amer. Math. Soc., Providence, RI, 2000. }
\pagestyle{myheadings}

\section{Introduction}\label{section-introduction}
\markboth{K. D. ELWORTHY AND XUE-MEI LI}{SPECIAL IT\^O MAPS AND ONE FORMS}

Let $M$ be a smooth compact Riemannian manifold. For a point $x_0$ of $M$
let $C_{x_0}M$ denote the space of continuous paths $\sigma:[0,T]\to M$
with $\sigma(0)=x_0$, for some  fixed $T>0$. Then $C_{x_0}M$ has a natural
$C^\infty$ Banach manifold structure, as observed by J. Eells, with
 tangent spaces $T_\sigma C_{x_0}M$ which can be identified with the 
spaces of continuous maps $v\colon [0,T]\to TM$ over $\sigma$ such that
$v(0)=0$, each $\sigma\in C_{x_0}M$. The Riemannian structure of $M$ 
induces a Finsler norm $\|\hskip4pt\|_\sigma^\infty$ on each
 $T_\sigma C_{x_0}M$
with $$\|v\|_\sigma^\infty=\sup\{|v(t)|_{\sigma(t)}: 0\le t\le T\}$$
so that $T_\sigma C_{x_0}M$, $\|-\|_\sigma^\infty$ is a Banach space.
 We can then form the dual spaces
 $\displaystyle{
T_\sigma^*C_{x_0}M={\Bbb L}\left(T_\sigma C_{x_0}M;{\Bbb R}\right)}$
to obtain the cotangent bundle $T^*C_{x_0}M$
whose sections are 1-forms on $C_{x_0}M$. To obtain $q$-vectors,
$0\le q<\infty$ take the exterior product $\wedge^q T_\sigma C_{x_0}M$
completed by the greatest cross norm 
\cite{Michor-LNM}
 so that the space of continuous linear maps 
${\Bbb L}\left(\wedge^q T_\sigma C_{x_0}M;{\Bbb R}\right)$
is naturally isomorphic to the space of alternating q-linear maps
$$\alpha: \hskip20pt
T_\sigma C_{x_0}M\times\dots\times  T_\sigma C_{x_0}M\longrightarrow \R$$
(and also to the corresponding completion $\wedge^q T_\sigma^*C_{x_0}M$).
Let $\Omega^q$ be the space $\Gamma\wedge^qT^*C_{x_0}M$ of sections of
the corresponding bundle. These are the q-forms. If $C^r\Omega^q$ refers
to the $C^r$ q-forms, $0\le r\le \infty$, then exterior differentiation
gives a map 
$$d:\hskip24pt  C^r\Omega^q \longrightarrow C^{r-1}\Omega^{q+1},
\hskip 50pt r\ge 1.$$
This is given by the formula: if $V^j$, $j=1$ to $q+1$, are $C^1$
vector fields, then for $\phi\in C^1\Omega^q$
\begin{equation}\label{Palais}
\begin{array}{ll}
&d\phi\left(V^1\wedge\dots \wedge V^{q+1}\right)\\
&=\sum_{i=1}^{q+1} (-1)^{i+1} L_{V^i}\left[\phi\left(V^1\wedge\dots 
\wedge \widehat{V^i} \wedge\dots \wedge V^{q+1}\right)\right]\\
&+\displaystyle{\sum_{1\le i<j\le q+1} (-1)^{i+j}  } 
 \phi\left( [V^i,V^j]\wedge  V^1\wedge\dots
  \widehat{V^i}\wedge \dots \widehat{V^j}\dots
 \wedge V^{q+1}\right)
\end{array}
\end{equation}
where $[V^i,V^j]$ is the Lie bracket and $ \widehat{V^j}$
means omission of the vector field $V^j$ , e.g. see
 \cite{Lang-manifold-book}.

\bigskip

For each $r\ge 1$ there are the deRham cohomology groups\\
 $H^q_{deRham(r)}\left(C_{x_0}M\right)$. If we  were using spaces of 
H\"older continuous paths, as in \cite{Bonic-Framphton-Tromba} we would have 
smooth partitions of unity and the deRham groups would be equal to the 
singular cohomology groups and so trivial for $q\ge 0$ since based path 
spaces are contractible.
An as yet unpublished result of C. J. Atkin carries
 this over to continuous paths, even though $C_{x_0}M$ does not admit smooth
 partitions of unity. In any case since our primary interest is in the 
differential analysis associated with the Brownian motion measure
 $\mu_{x_0}$ on $C_{x_0}M$, which could equally well be considered on
H\"older paths of any exponent smaller than a half, we could use H\"older
rather than continuous paths and it is really only for notational convenience
that we do not: the resulting manifold  would admit $C_\infty$
 partitions of unity \cite{Frampton-Tromba}.
 Independently of the existence  of partitions of unity:
contractibility need not imply triviality of the deRham cohomology
group when some restriction is put on the spaces of forms. For example if 
$f:\R\to \R$ is given by $f(x)=x$ then $df$ determines a non-trivial
class in the first bounded deRham group of $\R$. In finite dimensions
the $L^2$ cohomology of a cover $\tilde M$ of a compact manifold
$M$ gives important topological invariants of $M$ even when $\tilde M$
is contractible, eg see \cite{Atiyah76}; note also  
\cite{Bueler-Prokhorenkov}.

\bigskip

In finite dimensions the $L^2$ theory has especial significance because of
its relationship with Hodge theory and the associated geometric analysis.
In infinite dimensions L. Gross set the goal of obtaining an analogous
Hodge theory at the time of his pioneering work on infinite dimensional
potential theory \cite{GrossJFA1} in the late 60's. In his work he 
demonstrated the importance of the Cameron-Martin space $H$ in potential
analysis on Wiener space. In particular he showed that $H$-differentiability
was the natural concept in such analysis; a fact which became even
more apparent later, especially with the advent of Malliavin calculus.
For related analysis on infinite dimensional manifolds such as $C_{x_0}M$
the `admissible directions' for differentiability have to be subspaces
of the tangent spaces and the `Bismut tangent spaces', subspaces 
$H^1_\sigma$ of $T_\sigma C_{x_0}M$, revealed their importance in the
work of Jones-L\'eandre \cite{Jones-Leandre91}, and later in the integration
by parts theory of Driver \cite{Driver92} and subsequent surge of
 activity. They are defined by the parallel translation 
 $//_t(\sigma):T_{x_0}M\to T_{\sigma(t)}M$ of the Levi-Civita connection
and consist of  those $v\in T_\sigma C_{x_0}M$ such that 
$v_t=//_t(\sigma)h_t$ for $h_\cdot \in L_0^{2,1}\left(T_{x_0}M\right)$.
To have a satisfying $L^2$ theory of differential forms on $C_{x_0}M$
the obvious choice would be to consider `H-forms' i.e.
for 1-forms these would be $\phi$ with $\phi_\sigma\in (H^1_\sigma)^*$,
$\sigma\in C_{x_0}M$, and this agrees with the natural $H$-derivative
$df$ for $f:C_{x_0}M\to \R$. For $L^2$ q-forms the obvious  choice would be
$\phi$ with $\phi_\sigma\in\wedge^q(H_\sigma^1)^*$, using here the 
Hilbert space completion for the exterior product. An $L^2$-deRham theory
would come from the complex of spaces of $L^2$ sections
\begin{equation}\label{complex}
\dots\stackrel{\bar d}{\to}  L^2\Gamma\wedge^q(H^1_\sigma)^*
\stackrel{\bar d}{\to} L^2\Gamma\wedge^{q+1}(H^1_\sigma)^*
\stackrel{\bar d}{\to}\dots
\end{equation}
where $\bar d$ would be a closed operator obtained by closure from the usual 
exterior derivative (\ref{Palais}). From this would come the
 deRham-Hodge-Kodaira Laplacians $\bar d\bar d^*+\bar d^* \bar d$ and an
 associated Hodge decomposition. However the brackets $[V^i, V^j]$ of
sections of $H^1_\cdot$ are not in general  sections of $H_\cdot^1$,
and formula (\ref{Palais}) for $d$ does not make sense for $\phi_\sigma$
defined only on $\wedge^q H_\sigma^1$, each $\sigma$,
e.g. see \cite{Cruzeiro-Malliavin96}, \cite{Driver-Lie-bracket}.
 The project fails
at the stage of the definition of exterior differentiation. Ways to 
circumvent this problem were found, for paths and loop spaces, by L\'eandre
\cite{Leandre-cohomology96} 
 \cite{Leandre-homogeneous} \cite{Leandre-stochastic-cohomology}
 who gave analytical deRham groups and showed
 that they agree with the singular cohomology of the spaces. See also
\cite{Leandre-homology}. However these were
not $L^2$ cohomology theories and did not include a Hodge theory.
For flat Wiener space the problem with brackets does not exist
(the $H_\cdot^1$ bundle is integrable) and a full $L^2$ theory was
carried out by Shigekawa \cite{Shigekawa-Hodge}, including the
proof of triviality of the groups, see also \cite{Mitoma}.
For Wiener manifolds see \cite{Piech82}.
For paths on a compact Lie group $G$ with bi-invariant metric, and
 corresponding loop groups, there is an alternative, natural, 
 H-differentiability structure with the $H_\sigma$ modified by
using the flat left and right connection instead of the Levi-Civita
connection usually used. With this Fang\&Franchi \cite{Fang-Franchi97},
 \cite{ Fang-Franchi}, carried through a construction of the complex
(\ref{complex}) and obtained a Hodge decomposition for $L^2$ forms.
For a recent, and general, survey see \cite{Leandre-survey99}.

\bigskip

Our proposal is to replace the Hilbert spaces $\wedge^q H_\sigma^1$
in (\ref{Palais}) by other Hilbert spaces ${\mathcal H}_\sigma^q, q=2,3,\dots$,
continuously included in $\wedge^q T_\sigma C_{x_0}M$, though keeping the
exterior derivative a closure of the one defined by (\ref{Palais}). Here
we describe the situation for $q=2$ which enables us to construct an 
analogue  of the de Rham-Hodge-Kodaira Laplacian on $L^2$ sections
of the Bismut tangent spaces and a Hodge decomposition of the space of
$L^2$ 1-forms. At the time of writing this the situation for higher order 
forms is not so clear and the discussion of them, and some details of the
construction here, are left to a more comprehensive article.

\bigskip

The success of Fang and Franchi for path and loop groups was due to
a large extent to the fact that the It\^o map (i.e. solution map) of
the right or left invariant stochastic differential equations for
 Brownian motion on their  groups is particularly well behaved,
in particular its structure sends the Cameron-Martin space to the Bismut
type tangent space. The basis for the analysis here is the fact that
the It\^o map of gradient systems is almost as good: the `almost'
being made into precision by `filtering out redundant noise' 
\cite{Elworthy-Yor} \cite{Elworthy-LeJan-Li1}
\cite{Elworthy-LeJan-Li-book}.  
Indeed this is used for $q\ge 2$ to \underline{define} the spaces
${\mathcal H}^q_\sigma$, though it turns out that they depend only
on the Riemannian  structure of $M$, not on the embedding used to 
obtain the gradient stochastic differential equation. The good
properties of the It\^o maps have been used for analysis on path
and loop spaces, particularly by Aida, see \cite{Aida-Elworthy95}
\cite{Aida96}: Theorem  \ref{th:basic} consolidates these and should be of 
independent interest. The result that the It\^o map can be used to
continuously pull back elements of $L^2\Gamma({{\mathcal H}^1_\cdot}^*)$,
 i.e. $L^2$ H-forms to $L^2$ H-forms on Wiener space seems rather surprising.

\bigskip

We  should also mention the work done on `submanifolds' of Wiener space 
  and in particular on the submanifolds which give a model for the
based loop space of a Riemannian manifold. For this see 
\cite{Airault-VanBiesen-codimension} 
\cite{Kusuoka-cohomology} \cite{Kusuoka-forms}
\cite{VanBiesen-divergence}. For a detailed analysis of some analogous
properties of the stochastic development Ito map see 
\cite{X-D-Li-thesis}.

\bigskip

As usual in this subject all formulae have to be taken with the convention
 that equality only holds for all paths outside  some  set of measure zero.

\bigskip

\section{The It\^o map for gradient Brownian dynamical systems}
\label{section-Ito}

{\bf A.} Let $j: M\to \R^m$ be an isometric embedding. The existence
of such a $j$ is guaranteed by Nash's theorem. Let $X:M\times\R^m\to TM$
be the induced gradient system, so $X(x): \R^m\to T_xM$ is the orthogonal
projection, or equivalently $X(x)(e)$ is the gradient of 
$x \mapsto \langle j(x),e\rangle _{\R^m}$, for $e\in \R^m$. Take the
canonical Brownian motion $B_t(\omega)=\omega(t)$, $0\le t \le T$, for
$\omega \in \Omega=C_0(\R^m)$ with Wiener measure ${\Bbb P}$. The 
solutions to the stochastic differential equation

\begin{equation}\label{sde}
dx_t=X(x_t)\circ dB_t
\end{equation}
on $M$ are well known to be Brownian motions on $M$, \cite{Elworthy-book}
\cite{Elflour} \cite{Rogers-Williams-2}. Let  
 $\xi _t(x,\omega)$: $0\le  t\le T, x\in M, \omega \in \Omega$
denote its solution flow and

$${\mathcal  I}:C_0\left( {\Bbb R}^m\right) \rightarrow C_{x_0}(M), $$
${\mathcal  I}(\omega)_t=\xi_t(x_0, \omega)$ its It\^o map. Then $\I$ maps 
${\Bbb P}$ to the Brownian measure $\mu_{x_0}$ on $C_{x_0}M$.
The flow is $C^\infty$ in $x$ with random derivative
 $T_{x_0}\xi_t:T_{x_0}M\to T_{x_t}M$ at $x_0$. The It\^o map is smooth in
 the sense of Malliavin, as are all such It\^o maps, e.g. see
\cite{Ikeda-Watanabe} and \cite{Malliavin-book}, 
with H-derivatives continuous linear maps
$$T_\omega{\mathcal  I}: H\to T_{x_\cdot(\omega)}C_{x_0}M, \hskip 30pt
\hbox{almost all } \omega \in\Omega.$$
Here $x_\cdot(\omega)\ := \ \xi_\cdot(x_0,\omega)$ and $H$ is  the
 Cameron-Martin space $L_0^{2,1}(\R^m)$ of $C_0(\R^m)$.
From Bismut \cite{Bismut-Durham}  there is the formula
\begin{equation}\label{Ito-derivative}
T\I(h)_t=T_{x_0}\xi_t  \int_0^t (T_{x_0}\xi_s)^{-1}  X(x_s)(\dot h_s)ds
\end{equation}
for $h\in H$.

\bigskip

\underline{Remark:}
All the following results remain true when the gradient SDE (\ref{sde})
is replaced by an SDE of the form (3) with smooth coefficients whose
 LeJan-Watanabe connection, in the sense of \cite{Elworthy-LeJan-Li-book},
 is the Levi Civita connection of our Riemannian manifold (and consequently
whose solutions are Brownian motions on M). The only exception is the
 reference to the shape operator in \S\ref{section-decomposition}B below.
 The canonical stochastic differential equation for compact  Riemannian
 symmetric spaces gives a class of such stochastic differential equations,
 see \cite{Elworthy-LeJan-Li-book} section 1.4, with compact Lie groups
 giving specific examples.

\bigskip

{\bf B.} Formula (\ref{Ito-derivative}) is derived from the covariant
stochastic differential equation along $x_\cdot$ for $v_t=T\I(h)_t$ 
using the Levi-Civita connection

\begin{equation}\label{eq-derivative}
Dv_t=\nabla X_{v_t} \left(\circ dB_t\right)+X(x_t)(\dot{h}_t)dt,
\end{equation}
where ${D\over \partial t}=//{d\over dt}//^{-1}$ and $D$ is the
 corresponding stochastic differential. Clearly $v_\cdot$ does not lie in 
the Bismut tangent spaces in general. In fact $\nabla X$ is determined
by the shape operator of the embedding:
$$A: TM\times \nu M\to TM$$
where $\nu M$ is the normal bundle of $M$ in $\R^m$. We have
\begin{equation}\label{shape-op}
\nabla_vX(e)=A(v, K_xe), \hskip 55pt v\in T_xM
\end{equation}
where $K_x:\R^m \to \nu_xM$ is the orthogonal projection, $x\in M$.
We can think of $\nu_xM$ as $Ker(x)$. For $v_\cdot$ to be a Bismut
tangent vector for all $h$ would require $A\equiv 0$ so that $M$
would be isometric to an open set of $\R^n$.

\bigskip

Equation (\ref{eq-derivative}) has the It\^o form
\begin{equation}\label{eq-derivative-Ito-form}
 Dv_t=\nabla X_{v_t} \left(dB_t\right)-{1\over 2}Ric^{\#}v_t dt
+X(x_t)(\dot{h}_t)dt
\end{equation}
where $Ric^{\#}: TM\to TM$ corresponds to the Ricci tensor. From 
 (\ref{shape-op}) this  is driven only by the
 'redundant noise' in the kernel of $X(x_t)$. The technique of
 \cite{Elworthy-Yor} shows that if ${\mathcal  F}^{x_0}$ denotes the 
$\sigma$-algebra generated by $x_s: 0\le s \le T$ and 
\begin{equation}\label{filtered-derivative}
\bar v_t:=\overline{T{\mathcal  I}(h)}_t\ :=\ \E\{v_t\left|{\mathcal  F}^{x_0}\right.\}
\end{equation}
then
\begin{equation}\label{eq-derivtive-filtered}
{ D \over \partial t}\bar v_t 
    =-\frac 12 Ric^{\#}(\bar v_t)+X(x_t)(\dot h_t).
\end{equation}
This is described in greater generality in \cite{Elworthy-LeJan-Li-book}.
We will rewrite (\ref{eq-derivtive-filtered}) as
\begin{equation}\label{eq-derivtive-filtered-simple}
{{\Bbb D} \over \partial t}\bar v_t  =X(x_t)(\dot h_t)
\end{equation}
where 
 \begin{equation}
{{\Bbb D} \over \partial t}:
 ={D \over \partial t}+{1\over 2} Ric^{\#}(V_t).
\end{equation}
If $W_t:T_{x_0}M\to T_{x_t}M $ is the Dohrn-Guerra, or `damped' parallel
transport defined by 

$${{\Bbb D} \over \partial t}\left(W_t(v_0)\right) =0$$
then ${{\Bbb D} \over \partial t}=W_t{d\over dt}W_t^{-1}$. It appears to
be of basic importance e.g. see \cite{Meyer},  \cite{Nelson-book}, 
\cite{Leandre93}, \cite{Cruzeiro-Fang-95}, \cite{Norris-Brownian-sheet},
\cite{Malliavin-book},  \cite{Fang98},  \cite{Elworthy-LeJan-Li-book}. 
To take this into account we will change the interior product of the
 Bismut tangent spaces $H^1_\sigma$ and take 
\begin{equation}\label{inner-product}
<u^1,u^2>_\sigma \hskip 4pt :=      \int_0^T
 \left\langle {{\Bbb D} \over \partial s}u^1_s, \hskip4pt
 {{\Bbb D} \over \partial s} u^2_s\right\rangle ds.
\end{equation}
Let ${\mathcal  H}^1_\sigma$ denote $H_\sigma^1$ with this inner product:
it consists of those tangent vectors $v_\cdot$ above $\sigma_\cdot$ such 
that
  $\displaystyle{\int_0^T 
\left|{{\Bbb D}\over \partial t}v_t\right|^2_{\sigma(t)} \, dt <\infty}$.
Since $X(x_t)$ is surjective we see that
 $\displaystyle{h\mapsto \overline{T_\omega{\mathcal  I}}(h_\cdot) }$
maps $H$ onto ${\mathcal  H}^1_{x_\cdot(\omega)}$ for each $\omega$, and

\begin{equation}\label{Ito-isometry}
\left|\left| \overline {T\I(h)} \right|\right|_{x_\cdot(\omega)}=
\sqrt {\int_0^T \left| X(x_t(\omega)) \dot h_t \right|^2 dt}.
\end{equation}
Let $\overline{T\I}()_\sigma: H\to {\mathcal  H}^1_\sigma$ denote the map
\begin{equation}\label{filtered-derivative-notation}
h\mapsto \E\left\{T\I_\cdot(h_\cdot)\left|
 x_\cdot(\omega)=\sigma     \right.  \right\}
\end{equation}
defined for $\mu_{x_0}$ almost all $\sigma\in C_{x_0}M$.

\bigskip

When $h:C_0(\R^m)\to H$ gives an adapted process we see, e.g. from
 \cite{Elworthy-LeJan-Li-book} that
$$ \overline{TI(h)}_\sigma\ :=\ 
\E\left\{T\I_\cdot(h_\cdot)\left|x_\cdot=\sigma\right.  \right\}
=\overline{T\I}(\bar h_\cdot (\sigma))_\sigma$$
where $\bar h(\sigma)_t=\E\left\{h_t\left| x_\cdot(\omega)=\sigma
              \right.  \right\}$. 
For non-adapted $h:C_0(\R^m)\to H$ see Theorem \ref{th:basic} below.

\vskip 15pt

{\bf C.} Let $\phi$ be a $C^1$ 1-form on $C_{x_0}M$ which is bounded together
with $d\phi$, using the Finsler structure defined above. For example
$\phi$ could be cylindrical and $C^\infty$. Then there is the pull
back $\I^*(\phi):C_0(\R^m)\to H^*$,  the H-form on $C_0(\R^m)$  given by
\begin{equation}\label{pull-back-def}
\I^*(\phi)_\omega(h)=\phi\left(T_\omega\I(h)\right), 
\hskip24pt h\in H.
\end{equation}
Also  $\I^*(\phi)\in D(\bar d)$, for $\bar d$ the closure of the exterior
derivative on $H$-forms on Wiener space and 
\begin{equation}\label{dI=Id}
 \bar d\left(\I^*(\phi)\right) =\I^*\left (d\phi\right),
\end{equation}
for $\I^*$ the pull back
$$\I^*(\psi)_\omega(h^1\wedge h^2)=
    \psi\left(T_\omega\I(h^1)\wedge T_\omega\I(h^2)\right), $$
 when $\psi \in \Omega^2$ and  $h^1,h^2\in H$, 
so $\I^*(\psi):C_0(\R^m)\to (H\wedge H)^*$, c.f. \cite{Malliavin-book}
and \cite{Fang-Franchi}. This can be proved by approximation.

\bigskip

We will show that $\I^*(\phi)$ can be defined when $\phi$ is only an
H-form on $C_{x_0}M$ although the right hand side of (\ref{pull-back-def})
is not, classically, defined in this case.

\vskip 15pt

{\bf D.} For each element $h$ of a Hilbert space $H$, let $h^\#$ denote
the dual element in $H^*$, and conversely, so $(h^\#)^\#=h$.
For Hilbert spaces $H$, $H^\prime$ and measure spaces
 $(\Omega, \F,{\Bbb P})$,  $(\Omega^\prime, \F^\prime,{\Bbb P}^\prime)$,
 a linear map 
$$S: L^2\left(\Omega^\prime, \F^\prime, {\Bbb P}^\prime; {H^\prime}^*\right)
  \to L^2\left(\Omega, \F, {\Bbb P}; H^*\right)$$
will be said to be the \underline{co-joint} of a linear map
$$T: L^2\left(\Omega, \F, {\Bbb P}; H\right)
\to  L^2\left(\Omega^\prime, \F^\prime, {\Bbb P}^\prime; H^\prime\right)$$
if 
$$\displaystyle{ h\mapsto \left[S(h(\cdot)^{\#})(\cdot)\right]^{\#}}$$
is the usual Hilbert space adjoint $T^*$ of $T$, i.e. if
$$\int S(\phi)_\omega h(\omega) d{\Bbb P}(\omega)
=\int \phi_\omega\left(T(h)(\omega)\right)   d{\Bbb P}^\prime(\omega)$$
all $\phi\in L^2(\Omega^\prime;H^\prime)$ and $h\in L^2(\Omega; H)$.

A linear map $T$ of Hilbert spaces is a {\it Hilbert submersion} if $TT^*$
is the identity map.

\bigskip

We first note a preliminary result
\begin{proposition}
The map $\displaystyle{ \overline{T\I}:
 L^2\left(C_0(\R^m),\F^{x_0},{\Bbb P}|_{\F^{x_0}}; H\right)
 \to L^2\Gamma{\mathcal  H_\cdot}^1} $
given by  $\displaystyle{\overline{T\I}(h_\cdot)(\sigma)
=\overline{T\I}(h_\cdot(\sigma))(\sigma)}$ 
is 
\begin{equation}\label{formula-filtered-derivative}
h\mapsto W_\cdot\int_0^\cdot W_s^{-1} X(\sigma(s))\dot h(\sigma)_s ds
\end{equation}
and is a Hilbert submersion with inverse and adjoint given by
\begin{equation}\label{filtered-derivative-inverse}
v \mapsto \int_0^\cdot Y(x_s(\cdot)){{\Bbb D}\over \partial s} v_s ds
\end{equation}
for $Y(x): T_xM\to \R^m$ the adjoint of $X(x)$, $x\in M$. Its co-joint
can be written 
$\displaystyle{\overline{\I^*(-)}: L^2\Gamma({{\mathcal  H}^1_\cdot}^*)
\to L^2\left(C_0(\R^m), \F^{x_0};H^*\right)}$ in the sense that it agrees
with $\phi\mapsto \E\{\I^*(\phi)\left| \F^{x_0}\right. \}$
for $\phi$ a 1-form on $C_{x_0}M$.
\end{proposition}

{\it Proof.}
Note $h\mapsto X(\sigma_\cdot)(\dot h_\cdot)$ maps $h$ to the space
$L^2T_\sigma C_{x_0}M$ of $L^2$ `tangent vectors' to $C_{x_0}M$ at
$\sigma$, and as such is a Hilbert projection with inverse and adjoint
$u\mapsto \int_0^\cdot Y(\sigma_s)(u_s)ds$ since
\begin{eqnarray*}
\langle\int_0^\cdot Y(\sigma_s)u_sds, h_\cdot\rangle_H
&=&\int_0^T \langle Y(\sigma_s)u_s, \dot h_s\rangle_{\R^m} ds\\
&=&\int_0^T \langle u_s, X(\sigma_s)\dot h_s\rangle_{\sigma_s} ds\\
&=&\langle u_\cdot, X(\sigma_\cdot)\dot h_\cdot
  \rangle_{L^2 T_\sigma C_{x_0}M}.
\end{eqnarray*}
Also $\displaystyle{u\mapsto W_\cdot\int_0^\cdot W_s^{-1}u_s\  ds}$\ 
is an isometry of $L^2T_\sigma C_{x_0}M$ with ${\mathcal  H}^1_\sigma$ 
by definition, and its inverse is ${{\Bbb D}\over \partial s}$.
To check the co-joint:
if $h\in L^2(C_0(\R^m),\F^{x_0};H)$ and $\phi\in \Omega^1$  is bounded and
 continuous then
\begin{eqnarray*}
\int_{C_0(\R^m)} \overline{\I^*(\phi)}(h)\ d{\Bbb P}
&=& \int_{C_0(\R^m)} \I^*(\phi)(h)\  d{\Bbb P}\\
&=&\int_{C_0(\R^m)} \phi(T\I (h))\  d{\Bbb P}\\
&=&\int_{C_0(\R^m)} \phi\left( \E\left\{T\I (h)\ \ \left|\ \ \F^{x_0}
 \right.\right\} \right)\ 
 d{\Bbb P}\\
&=&\int_{C_{x_0}(M)} \phi\left(
    \overline{T\I}(h)\right) \ d\mu_{x_0}.\\
\end{eqnarray*}
\hfill Q.E.D.

\bigskip

Our basic result on the nice behaviour of our It\^o map is the following:
\begin{theorem}\label{th:basic}
The map $h\mapsto \E\left\{T\I(h)\left| \F^{x_0}\right.\right\}$
 determines a continuous linear map 
$$\overline{T\I(-)}: L^2\left(C_0(\R^m);H\right)\to
L^2\Gamma {\mathcal  H}^1,$$
which is surjective. The pull back map $\I^*$ on 1-forms extends to
a continuous linear map of H-forms:

$${\mathcal  I}^{*}:L^2\Gamma ({{\mathcal  H}^1_\cdot}^{*})\rightarrow L^2\left(
C_0\left( {\Bbb R}^m\right) ;H^{*}\right),  $$
which is the co-joint of $\overline{T\I(-)}$. It is injective with
closed range.
\end{theorem}
The proof of the continuity of $\overline{T\I(-)}$ is given in the next
section (\S\ref{section-decomposition}D).
 Its surjectivity follows from the previous proposition.
 That its co-joint agrees  with $\I^*$ on $\Omega^1$ comes from
the last few lines of the proof  of that proposition. From this we have
the existence of the claimed extension of $\I^*$ and
its continuity, injectivity and the fact that it has closed range.

\vskip 15pt

{\bf E.} Let
 $\bar d: Dom(\bar d)\subset L^2(C_0(\R^m);\R)\to  L^2(C_0(\R^m);H^*)$
 be the usual closure of the H-derivative, as in
Malliavin calculus. Let $\bar d:Dom(\bar d)\subset L^2(C_0(M);\R)\to 
 L^2\Gamma{{\mathcal  H}^1}^*$ be the closure of differentiation in
 ${\mathcal  H}^1_\cdot$   directions defined on smooth cylindrical
 functions (to make a concrete choice). The existence of this closure is
 assured and well known by Driver's integration by parts formula. We note
 the following consequence of Theorem \ref{th:basic}, although it is
 not needed in the following sections. In it we also use $\I$ to pull
 back functions on $C_{x_0}M$ by
  $\I^*(f)(\omega)=f(\I^*(\omega))=f(x_\cdot(\omega))$.

\begin{corollary}\label{co:Ito-d}
With $\I^*$ defined on H-forms by  Theorem \ref{th:basic} the compositions
$\I^*\bar d$ and $\bar d \I^*$ are closed, densely defined operators
on their domains in $L^2(C_{x_0}M;\R)$ into $L^2(C_0(\R^m);H^*)$ and
\begin{equation}\label{Ito-d}
\I^*\bar d\subset \bar d \I^*.
\end{equation}
\end{corollary}

{\it Proof.} 
Let Cyl denote the space of smooth cylindrical functions on $C_{x_0}M$.
 If $f\in Cyl$ it is standard that $\I^*(f)\in Dom(\bar d)$  so 
$Cyl\subset Dom(\I^*\bar d)\cap Dom(\bar d\I^*)$ since by  Theorem 
\ref{th:basic} $Dom(\I^*\bar d)=Dom(\bar d)$. Moreover
$$\I^*df=\bar d \I^*f$$
by the chain rule. Since $\I^*$ is continuous on functions $\bar d I^*$
 is closed and we have 
$$\overline{\I^*d|_{Cyl}}^c=\overline{\bar d\I^*|_{Cyl}}^c
  \subset \bar d\I^*$$
where $\overline{\\}^c$ denotes closure. Indeed $\I^*\bar d$ is closed
since $\I^*$  is continuous with closed range on H-forms
by the theorem  so the closure $\overline{\I^*d|_{Cyl}}^c$  exists and
 is a restriction of $\I^*\bar d$. The result follows by showing this
 restriction is in fact equality: For this suppose $f\in Dom(\I^*\bar d)$.
Then $f\in Dom(\bar d)$ so there exists $f_n\in Cyl$ with $f_n\to f$ in
$L^2$ and $df_n\to \bar d f$. By continuity of $\I^*$ we have
 $\I^*(df_n)\to \I^*(\bar df)$ so that 
$f\in Dom(\overline{\I^*d|_{Cyl}}^c)$. \hfill Q.E.D.

\bigskip

Taking co-joints, and defining $-div$ to be the co-joint of $\bar d$
on $C_{x_0}M$ and $C_0(\R^m)$, we have the following corollary.
It formalises some of the arguments in \cite{EL-LIibp} 
\cite{Elworthy-LeJan-Li-book}.

\begin{corollary}
The composition $div\, \overline{T\I(-)}$ and 
$\ \E\,\{div\left|\F^{x_0}\right.\}$ are closed densely defined operators
on $L^2(C_0(\R^m); H)$ into $L^2(C_{x_0}M;\R)$. Moreover
\begin{equation}\label{divergence-comp}
\E\{\,div\, \left|\,\F^{x_0}\right.\}\subset div\,\overline{T\I(-)}.
\end{equation}
\end{corollary}

Also by  a comparison result of H\"ormander, see \cite{Yosida-book},Thm2,
\S 6 of Chapter II, p79, (\ref{Ito-d})
implies there exists a constant $C$ such that 
\begin{equation}\label{comparison}
\int_{C_0(\R^m)} |\bar d \I^*(f)|^2_{H^*} d{\Bbb P}
\le C\left(\int_{C_{0}(\R^m)} |\I^* (\bar df)|^2_{H^*} d{\Bbb P}
+\|f\|_{L^2}^2\right)
\end{equation}
for all $f\in Dom(\I^*\bar d)=Dom(\bar d)$ in $L^2(C_{x_0}M;\R)$.

Let $L^{2,1}$ denote the domain of the relevant $\bar d$ with its graph 
norm  $$\|f\|_{{2,1}} =\sqrt{\|\bar df\|^2_{L^2}+\|f\|^2_{L^2}},$$
i.e. the usual Dirichlet space. The boundedness of $\I^*$ in the next
 corollary was proved directly for cylindrical functions (and hence for 
all $f\in L^{2,1}$) by Aida and Elworthy as the main step in their
proof of the logarithmic Sobolev inequality on path spaces.

\begin{corollary} c.f. \cite{Aida-Elworthy95}
The pull back determines a continuous linear map
$\I^*: L^{2,1}(C_{x_0}M)\to L^{2,1}(C_0(\R^m))$. It is injective with
closed range.
\end{corollary}

{\it Proof.} Continuity and existence is immediate from (\ref{comparison}) and
the continuity of $\I^*$ in Theorem \ref{th:basic}. Injectivity is clear.
To show the range is closed suppose $\{f_n\}_{n=1}^\infty$ is a sequence
in $L^{2,1}(C_{x_0}M)$ with $\I^*(f_n)\to g$, in $L^{2,1}(C_0(\R^m))$
some $g$.

 Then $g$ is $\F^{x_0}$ measurable so $g=\I^*(\bar g)$ some 
$\bar g\in L^{2}(C_{x_0}M)$. Moreover $f_n\to \bar g$ in 
$L^2(C_{x_0}M)$. We have $\bar d\I^*(f_n)\to \bar d g$. By
 (\ref{Ito-d}) $\bar d \I^*(f_n)=\I^*(\bar d f_n)$
 since $f_n\in Dom(\bar d)$. From this we see $\bar d f_n$ converges in
 $L^2\Gamma {\mathcal  H}_\cdot$ because $\I^*$ has closed range in 
$L^2(C_0(\R^m);H)$ by Theorem \ref{th:basic}. This shows
 $\bar g\in Dom(\bar d)$ as required.

\bigskip

\underline{Remark:}
We chose the basic domain of $d$ on functions of $C_{x_0}M$ to be
smooth cylindrical functions but the results above would hold
equally well with other choices e.g. the space $BC^1$ of functions
$F: C_{x_0}M\to \R$ which are $C^1$ and have $dF$ bounded on $C_{x_0}M$
using the Finsler structure of $C_{x_0}M$. The crucial conditions
needed for $Dom(d)$ are that it is dense in $L^2(C_{x_0}M)$ and that
$F\in Dom(d)$ implies $\I^*(F)\in Dom(\bar d)$ on Wiener space. At present
it seems unknown as to whether different choices  give the same closure
 $\bar d$. This is  essentially equivalent to knowing that the closed
 subspace  $\I^*(L^{2,1}(C_{x_0}M))$ of $L^{2,1}(C_0(\R^m))$ is
 independent of the choice of $Dom(d)$. The obvious guess would be that
 it is and consists of the $\F^{x_0}$-measurable elements of
 $L^{2,1}(C_0(\R^m))$, but we do not pursue that here.

\section{Decomposition of noise; proof of Theorem}
\label{section-decomposition}

{\bf A.} 
Let $\R^m$ denote the trivial bundle $M\times \R^m\to M$.
It has the subbundle $ker X$ and its orthogonal complement $Ker X^\perp$.
The projection onto these bundles induce connections on them, 
\cite{Elworthy-LeJan-Li1} and these combine to give parallel translations
$$\tilde{//_t}(\sigma):\R^m \to \R^m, \hskip 28pt  0\le t\le T$$
along almost all $\sigma \in C_{x_0}M$, which map $Ker X(x_0)$ to
$Ker X(\sigma(t))$ and preserve the inner product of $\R^m$. From
\cite{Elworthy-Yor} (extended to more general stochastic differential
equations in \cite{Elworthy-LeJan-Li-book}), there is a Brownian motion
$\{\tilde B_t: 0\le t\le T\}$ on $Ker X(x_0)^\perp$ and one,
 $\{\beta_t: 0\le t\le T\}$ on $Ker X(x_0)$ with the property that
\begin{enumerate}
\item
$\tilde B_\cdot$ and $\beta_\cdot$ are independent;
\item \label{ii} $\sigma(\tilde B_s: 0\le s\le t)=\sigma(x_s: 0\le s\le t)$,
 $0\le t\le T$ and in particular $\sigma(\tilde B_s: 0\le s\le T)=\F^{x_0}$;
\item $dB_t=\tilde{//}_td\tilde B_t+\tilde{//}_td\beta_t$.
\end{enumerate}

Let $\displaystyle{L^2\left(\F^{x_0};\R\right)}$, 
$\displaystyle{L^2\left(\beta;\R\right)}$ etc. denote the Hilbert 
subspaces of  $\displaystyle{L^2\left(C_0(\R^m);\R\right)}$ etc. 
consisting  of elements measurable with respect to $\F^{x_0}$
or $\sigma\{\beta_s: 0\le s\le T\}$. By (\ref{ii}) above
$\displaystyle{L^2\left(\F^{x_0};\R\right)= L^2\left(\tilde B;\R\right)}$.
As before we can identify  $\displaystyle{L^2\left(\F^{x_0};\R\right)}$
with $\displaystyle{L^2\left(C_{x_0}M;\R\right)}$.

\begin{lemma}\label{lemma:isomorphism1}
The map $f\otimes g\mapsto f(\cdot)g(\cdot)$
determines an isometric isomorphism
$$L^2\left(\F^{x_0};\R\right)\hat\otimes L^2\left(\beta;\R\right)
\to L^2\left(C_0(\R^m);\R\right)$$
where $\hat\otimes$ denotes the usual Hilbert space completion.
\end{lemma}
{\it Proof.} This is immediate from the independence of 
 $\tilde B_\cdot$ and $\beta$, the fact that 
$\tilde B_\cdot\times \beta: C_0\left(\R^m\right) \to C_0(Ker X(x_0)^\perp)
\times C_0(Ker X(x_0))$
generates $\F$, and the well known tensor product decomposition of $L^2$
of a product space.
\hfill Q.E.D.
\bigskip

{\bf B.} Let us recall the representation theorem for Hilbert space
valued Wiener functionals:

\begin{lemma}[$L^2$ representation theorem]
Let $\left\{ \beta _t, 0\le t\le T\right\} $ be an $\displaystyle{m-n}$
 dimensional  Brownian motion and $H$ a separable Hilbert space.
 Let $K$ be the Hilbert space of $\beta$-predictable
 ${\Bbb L}\left( {\Bbb R}^{m-n};H\right)$ valued process with

$$\left\| \alpha \right\| _K=\sqrt
{\int_0^T{\Bbb E\,}\left\| \alpha _r\left( \cdot\right) 
 \right\| _{{\Bbb L}\left( {\Bbb R}^{m-n},H\right)}^2\,dr}<\infty
$$
(using the Hilbert-Schmidt norm on
 ${\Bbb L}\left( {\Bbb R}^{m-n},H\right)$).  Then the  map

\begin{eqnarray*}
H\times K &\to &L^2(\beta; H)  \\
(h,\alpha)  &\mapsto &\int_0^T  \alpha_r(\cdot) (d\beta_r) +h
\end{eqnarray*}
is an isometric isomorphism. 
\end{lemma}

{\it Proof.} That it preserves the norm is a basic property of the
 Hilbert space valued It\^o integral. To see that it is surjective
 just observe that the image of the set of $\alpha\in K$ of the form
$$\alpha_s(\omega) = g(\omega) h_s$$
for  $g\in L^2(\beta;\R)$, $h\in H$  is total in  $L^2(\beta;H)$ by 
the usual representation theorem for real valued functionals and the
 isometry of $L^2(\beta;\R)\hat\otimes H$ with $L^2(\beta;H)$.
\hfill Q.E.D.

\bigskip

\begin{lemma}\label{lemma:isometry}
The map 
\begin{eqnarray*}
L^2\left( \beta;\R \right) \hat{\otimes}H &\longrightarrow &L^2\Gamma \left(
{\mathcal  H}_\cdot^1\right)  \\
g\otimes h &\mapsto &\overline{T{\mathcal  I}(gh)}
\end{eqnarray*}
is continuous linear.
\end{lemma}

\underline{Remark:}  Note that 
 $$\overline{T{\mathcal  I}(gh)}_\sigma
=\E\left\{T\I(gh)\ \left|\  x_\cdot=\sigma\right.\right\}
=\E\left\{gT\I(h)\ \left|\  x_\cdot=\sigma\right.\right\}.$$

{\it Proof.}  
By the representation theorem, a typical element $u$ of
 $L^2\left(\beta;H\right)$ has the form

$$u_t= h_t +\int_0^T\alpha _r(t) \left(d\beta _r\right), \hskip 24pt
0\le t\le T$$
for $h\in H=L_0^{2,1}(\R^m)$ and $\alpha \in K$; writing $\alpha_r(t)(e)$
for $\alpha_r(e)_t$, $e\in Ker X(x_0)^\perp \approxeq \R^{m-n}$.
 Now by equation (\ref{Ito-derivative})

\begin{eqnarray*}
&&\overline{ T{\mathcal  I}_t\left( \int_0^T\left\langle \alpha _r(\cdot),
 d\beta _r\right\rangle \right) } \\
&=&{\Bbb E}\left\{\ \  {\Bbb E}\left\{ T{\mathcal  I}_t\left( 
\int_0^T\alpha_r(\cdot) ( d\beta _r) \right) \ \ 
\left|\ \  {\mathcal  F}^{x_0}\vee {\mathcal  F}_t\right.\ \  \left|\ \ 
  {\mathcal  F}^{x_0}\right.
\right\} \right\} \\
&=&{\Bbb E}\ \left\{ T\xi _t\int_0^tT\xi _s^{-1}X(x_s) \left(
\int_0^t \dot\alpha_r(\cdot) (d\beta _r)\right)_s ds\ \ 
\left|\ \  {\mathcal  F}^{x_0}\right. \right\},
\end{eqnarray*}
where $\dot \alpha_r(s)$ means the derivative with respect to $s$.
Set 
\begin{eqnarray*}
u_t&=& T\xi_t \int_0^t (T\xi_s)^{-1} X(x_s) \left( \int_0^t 
\dot \alpha_r(\cdot) (d\beta_r) \right)_s ds\\
&=& T\I_t\left(\int_0^t \alpha_r(\cdot)d\beta_r\right)
\end{eqnarray*}
Then, by equation (\ref{eq-derivative-Ito-form}), writing
$\nabla_v X(e)=\nabla X(v)(e)=\nabla X^e(v)$ 
we see that $u_t$ has covariant It\^o differential given by
\begin{eqnarray*}
Du_t &=&\nabla X(u_t) dB_t- {1\over 2}Ric^{\#}(u_t)dt+
  X(x_t)\left( \int_0^t\dot\alpha_r(\cdot)(d\beta_r)\right)_t dt\\
&&+   T\I_t  \left(\alpha_t(\cdot) (d\beta_t)\right)
 +{1\over 2} \sum_{i=1}^{m-n} \nabla X^{e^i}\left( T\I_t \left(
      \alpha_t(\cdot)(e^i)    \right)\right) dt \\
\end{eqnarray*}
where $e^1,e^2,\dots,e^{m-n}$ is an orthonormal basis of
 $Ker X(x_0)^\perp$.

By properties (1) and (2) of $\beta$ and $\tilde B$ given in
\S\ref{section-decomposition}A above we can argue as in \cite{Elworthy-Yor}
\cite{Elworthy-LeJan-Li-book} to see that
$\displaystyle{\bar u_t\ :=\E\{u_t\ \left|\ \F^{x_0}\right. \}}$
satisfies 
\[
D{\bar u_t}=-\frac 12Ric^{\#}\left({\bar u_t}\right) dt+\frac
12\sum_{i=1}^{m-n}\nabla X^{e^i}\left( \overline{T{\mathcal  I}_t\left( \alpha
_t\left( \cdot \right) \left( e^i\right) \right) }\right) dt. 
\]
Thus 
\begin{eqnarray*}
 4\ \E\left| \bar u_\cdot\right|_{{\mathcal  H}^1}^2 
 &=&4\ \E \int_0^T \left|{\D \bar u_t \over \partial t}\right|^2\ dt\\
 &=&{\Bbb E}\int_0^T  \left| \sum_{i=1}^{m-n}
 \nabla X^{e^i}\left( {\Bbb E}\left\{ T{\mathcal  I} _r\left(
  \alpha _r\left( \cdot \right) \left( e^i\right) \right)\ \ 
\left|\ \  {\mathcal  F}^{x_0}\right. \right\} \right) \right| ^2dr \\
&\le &\int_0^T{\Bbb E}\left[ {\Bbb E}\left\{ \sum_{i=1}^{m-n}\left|\nabla
X^{e^i}\left( T{\mathcal  I}_r\left( \alpha _r\left( \cdot \right)
\left( e^i\right) \right) \right)\right| 
 \ \ \left|\ \ {\mathcal  F}^{x_0}\right. \right\} \right]^2\ dr
 \\
&\le &const\ \int_0^T \E\left[ \sum_{i=1}^{m-n}\E\left\{ \left\|
\nabla X^{e^i}\left( T{\mathcal  I}_r(-)\right) \right\|^2
	_{{\Bbb L}(H;T_{x_r}M)}\ \ 
   \left|\ \  {\mathcal  F}^{x_0}\right. \right\}
	\right. \cdot\\
&&\hskip 10pt\left.  
   \E\left\{ \left| \alpha _r(\cdot)( e^i) \right|^2_H\ \ 
   \left|\ \  {\mathcal  F}^{x_0}\right. \right\} \right]\  dr
\\
&\le &const\sup_{0\le r \le T}\E\left(||T\I_r||^2_{{\Bbb L}(H;T_{x_r}M)}
 	\right)\ \ 
\int_0^T \sum_{1}^{m-n} \E|\alpha_r(\cdot)(e^i)|^2_H dr\\
&=&const\left\| \int_0^T\alpha _r\left( \cdot \right) d\beta _r\right\|
_{L^2}^2,
\end{eqnarray*}
since $T\I: H\to TC_{x_0}M$  is well known to be in $L^2$, using the 
Finsler metric of $C_{x_0}M$, as can be seen from formula
 (\ref{Ito-derivative}) and standard estimates for
 $|T\xi_t|$, $|T\xi_t|^{-1}$, e.g. \cite{Bismut-aleatoire}, \cite{Kifer88},
 or use \cite{Elworthy-book}.
\hfill Q.E.D.

\bigskip

{\bf D.} \underline{Completion of Proof of Theorem \ref{th:basic}}

$$\overline{T{\mathcal  I}()}: L^2\left( C_0\left( {\Bbb R}^m\right);H\right)
  \longrightarrow L^2\Gamma \left( {\mathcal  H}^1\right) 
$$
is continuous.

{\it Proof.} Consider  the sequence of maps

\begin{eqnarray*}
L^2 (C_0(\R^m); H) &\approxeq& L^2( C_0(\R^m);\R)\, \hat\otimes\,H\\
 &&\approxeq L^2(\F^{x_0};\R)\,\hat\otimes\, L^2(\beta;\R)
     \, \hat\otimes\,H \\
&&\stackrel{1\otimes \overline{TI(-)}}{\vector(1,0){70}}\ \ 
 L^2(\F^{x_0};\R)\,\hat\otimes\,L^2\Gamma({\mathcal  H}_\cdot^1) \\
&&\stackrel{f\otimes v \mapsto fv}{\vector(1,0){70}}\ \ 
 L^2\Gamma({\mathcal  H}_\cdot^1) 
\end{eqnarray*}
where the first is the natural isomorphism and the second the isomorphism
given by Lemma \ref{lemma:isomorphism1}. The maps factorise our map
 $\overline{T\I(-)}$ and are all continuous,
 using Lemma \ref{lemma:isometry}.
 \hfill Q.E.D.

\section{Exterior differentiation and the space  ${\mathcal  H}^2_\sigma$}

\label{section-exterior}

{\bf A.} From $T_\omega\I:H\to T_{x_\cdot(\omega)}C_{x_0}M$ we can
form the linear map of 2-vectors
$$\wedge^2(T_\omega\I):\wedge^2H\to \wedge^2 T_{x_\cdot(\omega)}C_{x_0}M$$
determined by
$$h^1\wedge h^2\mapsto T_\omega\I(h^1)\wedge T_\omega\I(h^1).$$

Here as always, we use the usual Hilbert space cross norm on $H\otimes H$
and $\wedge^2 H$ is the corresponding Hilbert space completion, while
for $T_\sigma C_{x_0}M\otimes T_\sigma C_{x_0}M$ we use the greatest
cross norm (in order to fit in with the usual definition of differential
forms as alternating continuous multilinear maps). To see that this map
exists and is continuous using these completions we can use
 \cite{Carmona-Chevet79}, or directly use the characterisation of
$\wedge^2H$ as a subspace of the functions
 $\underline h: [0,T]\times [0,T]\to \R^m\otimes\R^m$ such that 
$$\underline h(s,t)=\int_0^s\int_0^tk(s_1,t_1)ds_1dt_1$$
for some $k\in L^2\left([0,T]\times[0,T]; \R^m\otimes\R^m\right)$,
 and the characterisation of
 $\wedge^2T_\sigma C_{x_0}M$ as a space of continuous functions 
$V$ into $TM\times TM$ with

\noindent
$\displaystyle{
 (s,t)\mapsto V_{(s,t)}\in T_{\sigma(s)}M\otimes T_{\sigma(t)}M}$.

\bigskip

For $\underline h \in \wedge^2 H$ we can form
$$\overline{\wedge^2(T\I)}(\underline h)_\sigma \ :=\
\E\ \{\wedge^2(T\I)(\underline h)\  \left| \ x_\cdot(\omega)=\sigma \right.\}
\in \wedge^2(T_\sigma C_{x_0}M)$$
almost all $\sigma\in C_{x_0}M$. This gives a continuous linear map
 $$\overline{\wedge^2(T\I)}_\sigma:\ =\ \overline{\wedge^2(T\I)}()_\sigma:
    \  \wedge^2 H \   \to \wedge^2(T_\sigma C_{x_0}M).$$
Let ${\mathcal  H}^2_\sigma$ be its image with induced Hilbert space  structure
( i.e. determined by its linear bijection with  
  $\displaystyle{\wedge^2 H /  (Ker \overline{\wedge^2(T\I)}_\sigma )}$\ ).
 We quote the following without giving its proof here. For the sequel the
important point is that the spaces ${\mathcal  H}^2_\sigma$ are determined only
by the Riemannian structure of $M$ and are independent of the embedding
used to obtain $\I$. Note also that in general they will be distinct
from $\wedge^2{\mathcal  H}^1_\sigma$.

\begin{theorem}\cite{Elworthy-Li-Hodge2}

The space ${\mathcal  H}^2_\sigma$ consists of elements of
$\wedge^2T_\sigma C_{x_0}M$ of the form $V+Q(V)$ where
 $V\in \wedge^2{\mathcal  H}^1_\sigma$ and $Q: \wedge^2{\mathcal  H}_\sigma^1\to 
 \wedge^2_\sigma C_{x_0}M$ is the continuous linear map determined
by 
$$Q(V)_{(s,t)}=(1\otimes W^s_t)\  W_s^{(2)}
\int_0^s (W_r^{(2)})^{-1}  \ {\mathcal  R} (V_{(r,r)})\ dr,
\hskip 24pt 0\le s\le t\le T$$
where
(i) $W_s^t=W_t\ W_s^{-1}:T_{\sigma(s)}M\to T_{\sigma(t)}M$ for $W_t$
as in \S\ref{section-introduction}B,

(ii) $W_t^{(2)}:\wedge^2  T_{x_0}M\to \wedge^2  T_{\sigma(t)}M $ is the 
the damped translation of 2-vectors on $M$ given by

\begin{eqnarray*}
{D\over \partial t}W_t^{(2)}(u)&=&-{1\over 2}{\mathcal  R}^{(2)}_{\sigma(t)}
\left(W_t^{(2)}(u)\right)\\
W_0^{(2)}(u)&=&u,\hskip 44 pt u\in\wedge^2 T_{x_0}M,
\end{eqnarray*}
for  ${\mathcal  R}_{\sigma(t)}^{(2)}$ the Weitzenb\"ock curvature on 2-vectors
e.g. see  \cite{Ikeda-Watanabe}, \cite{Elflour},\cite{Elworthy-Yor},
 \cite{Elworthy-LeJan-Li-book}, \cite{Rosenbergbook} (where it is called
`the curvature endomorphism') and

(iii) ${\mathcal  R}: \wedge^2TM\to\wedge^2TM$ denotes the curvature operator.
\end{theorem}

In this theorem we have used the identification of $\wedge^2 T_\sigma C_{x_0}M$
with elements $V_{(s,t)}\in T_{\sigma(s)}M\otimes T_{\sigma(t)}$,
continuous in $(s,t)$, and with the natural symmetry property. Thus 
$Q(V)$ is determined by the values $Q(V)_{(s,t)}$ for $s<t$.

\bigskip

{\bf B.} By restriction, 2-forms on $C_{x_0}M$ i.e. sections of the dual 
bundle to $\wedge^2TC_{x_0}M$ can be considered as sections of 
$({\mathcal  H}^2_\cdot)^*$. Let $Dom(d)\subset L^2\Gamma({\mathcal  H}^1_\cdot)^*$
be a dense linear subspace consisting of differential 1-forms $\phi$
on $C_{x_0}M$ restricted to  the ${\mathcal  H}^1_\sigma$ such that the
exterior differential $d\phi$,  defined by (\ref{Palais}), restricts
to give an $L^2$ section of $({\mathcal  H}^2_\cdot)^*$. For example 
Dom(d) could consist of $C^\infty$ cylindrical forms on $C_{x_0}M$
or $C^1$ forms which are bounded together with $d\phi$, using the Finsler
 norms. This gives a densely defined operator
$$d\ :\ Dom(d)\subset L^2\Gamma({{\mathcal  H}^1_\cdot}^*)\to L^2\Gamma 
  ({{\mathcal  H}^2_\cdot}^*).$$
We will also assume that each $\phi\in Dom(d)$ satisfies
 $\I^*(\phi)\in Dom(\bar d)$ and $\I^*(d\phi)=\bar d\I^*(\phi)$,
where $\I^*(d\phi):=d\phi(\wedge^2(T\I)(-))$ and $\bar d$ refers to
the closure of exterior differentiation on Wiener space,
$$ \bar d\ :\ 
Dom(\bar d)\subset L^2(C_0(\R^m);H^*)\to L^2(C_0(\R^m);\wedge^2H^*)$$
as defined in \cite{Shigekawa-Hodge} or \cite{Malliavin-book}.
This condition is easily seen to be satisfied by the two examples of
 $Dom(d)$ just mentioned, for example by approximation of $\I$.

\begin{theorem}
The operator $d$ from $Dom(d)$ in
 $L^2\Gamma({{\mathcal  H}_\cdot^1}^*) \to L^2\Gamma({{\mathcal  H}^2_\cdot}^*) $
is closable with closure a densely defined operator $\bar d$,
$$\bar d: Dom(\bar d)\subset L^2\Gamma({{\mathcal  H}_\cdot^1}^*)
  \to L^2\Gamma({{\mathcal  H}_\cdot^2}^*) $$
\end{theorem}

{\it Proof.}
Suppose $\{\phi_n\}_{n=1}^\infty$ is a sequence in Dom(d) converging
 in $ L^2\Gamma({{\mathcal  H}_\cdot^1}^*)$ to $0$ with $d\phi_n \to \psi$ in
$ L^2\Gamma({{\mathcal  H}^2_\cdot}^*)$  for some $\psi$. It suffices to show
$\psi=0$. As usual the proof uses integration by parts;
 the following method is derived from one used for scalars as in 
\cite{Albeverio-Brasche-Roeckner}, \cite{Eberle}. Let $\lambda:C_{x_0}M\to\R$
be $C^\infty$ and cylindrical and let $\underline h\in \wedge^2H$.
Then 
\begin{eqnarray*}
&&\int_{C_{x_0}M} \lambda(\sigma)\psi\left(
   \overline{\wedge^2(T\I)}(\underline h)(\sigma)   \right)   \ d\mu_{x_0}\\
&=&\lim_{n\to \infty}\int_{C_{x_0}M} \lambda(\sigma) d\phi_n\left(
  \overline{\wedge^2(T\I)}(\underline h)(\sigma)     \right)  
         \ d\mu_{x_0}\\
&=&\lim_{n\to \infty}\int_{C_{x_0}M}  d\phi_n\left(
          \lambda(\sigma)  \overline{\wedge^2(T\I)}(\underline h)(\sigma) 
       \right)   \ d\mu_{x_0}\\
&=&\lim_{n\to \infty}\int_{C_0(\R^m)}  (d\phi_n)\left(
          \lambda(x_\cdot(\omega)) 
   \wedge^2(T_\omega\I)(\underline h)
       \right)   \ d\ {\Bbb P}(\omega)\\
&=&\lim_{n\to \infty}\int_{C_0(\R^m)} \I^* (d\phi_n)\left(
          \lambda(x_\cdot(\omega)) \underline h
       \right)   \ d\ {\Bbb P}(\omega)\\
&=&\lim_{n\to \infty}\int_{C_0(\R^m)} \bar d (\I^*\phi_n)\left(
          \lambda(x_\cdot(\omega)) \underline h
       \right)   \ d\ {\Bbb P}(\omega)\\
&=&\lim_{n\to \infty}\int_{C_0(\R^m)}  (\I^*\phi_n)\left(
         div \lambda(x_\cdot(\omega)) \underline h
       \right)   \ d\ {\Bbb P}(\omega)\\
&=&0
\end{eqnarray*}
by the continuity of $\I^*$ in Theorem \ref{th:basic}. Here we used
 the property that
 $ \I^* (d\phi_n)= \bar d (\I^*\phi_n)$ for $\phi_n\in Dom(d)$, and
have let 

$$-div: Dom(div)\subset L^2(C_0(\R^m);\wedge^2H)\to L^2(C_0(\R^m);H)$$
denote the co-joint of $\bar d$. From \cite{Shigekawa-Hodge},
$\I^*(\lambda)\underline h\in Dom(div)$, as can be seen explicitly
in this simple situation.

Since the smooth cylindrical functions are dense in $L^2$ the above shows
 that
$\psi(\overline{\wedge^2(T\I)}(\underline h)_\sigma)=0$
 for almost all $\sigma$. Since $\underline h$ was arbitrary,
  $\overline{\wedge^2(T\I)}(-)_\sigma$ maps onto ${\mathcal  H}^2_\sigma$ by
  definition, and $H$ is separable,  this implies $\psi=0$ a.s. as required.
\hfill Q.E.D.

\bigskip

Let $\bar d\equiv \bar d^1$ be the closure of $d$,
 using the previous theorem, and let $d^*=(\bar d^1)^*$ be its adjoint.

\bigskip

From now on we shall assume that $Dom(d)$ was chosen so that it contains
smooth cylindrical one forms. Since the basic domain of $d$ on functions
was taken to be smooth cylindrical functions this implies that $d$ maps
$Dom(d)$ to $Dom(d)$. This property making the following crucial fact
almost immediate:

\begin{proposition}\label{pr:d2=0}
If $f:C_{x_0}M\to \R$ is in $Dom(\bar d)$ then $\bar d f\in Dom(\bar d^1)$
and $\bar d^1\bar df=0$.
\end{proposition}

{\it Proof.}
Let $\{f_n\}_{n=1}^\infty$ be a sequence in $D(d)$ converging in $L^2$
to $f$ with $df_n\to \bar d f$ in $L^2\Gamma{\mathcal  H}^1_\cdot$.
 Then $d(df_n)=0$, since $d^2=0$, each $n$, so $\bar d(\bar df)=0$.
\hfill Q.E.D.

\bigskip

The proposition enables us to define the first $L^2$ cohomology
group $L^2H^1(C_{x_0}M)$ by 
$$L^2H^1(C_{x_0}M)={Ker\ \bar d^1\over Image\ \bar d}.$$

\bigskip

\underline{Remark:}
On functions $\bar d$ has closed range by the existence of a spectral
gap for $d^*d$ on functions, proved by Fang\cite{Fang-gap} and the functional
 analytical argument of H. Donnelly \cite{Donnelly-spectrum}
 Proposition 6.2.

\section{A Laplacian on 1-forms and Hodge decomposition Theorem}
\label{section-Laplacian}

{\bf A.} For completeness we go through the formal argument which gives a
  self-adjoint `Laplacian' on 1-forms and a Kodaira-Hodge decomposition. 
Define an operator $\tilde d$ on
 $L^2(C_{x_0}M;\R)
\oplus L^2\Gamma({{\mathcal  H}^1_\cdot}^*) \oplus L^2\Gamma({{\mathcal  H}^2_\cdot}^*)$ to itself by

$$Dom(\tilde d)
=Dom(\bar d)\oplus Dom(\bar d^1)\oplus L^2\Gamma({{\mathcal  H}^2_\cdot}^*)$$
and $\tilde d(f,\theta,\phi)=(0,df,\bar d\theta)$
for $(f,\theta,\phi)\in Dom(\tilde d)$.
\bigskip

The operator $\tilde d+\tilde d^*$ has domain  the intersection
of the two domains, i.e. 
 $$Dom(\tilde d+\tilde d^*)=Dom(\bar d)\oplus (Dom (\bar d^1)
 \cap Dom(d^*)\,) \oplus Dom ({\bar {d^1}}^*).$$
From Driver's integration by parts formula it is known that $Dom(d^*)$
contains all smooth cylindrical 1-forms, as does $Dom (\bar d^2)$
by assumption. Thus $\tilde d+\tilde d^*$ has dense domain.
It is clearly symmetric. Furthermore, using Proposition \ref{pr:d2=0},
the consequent orthogonality of $Image(\bar d)$ and $Image(d^1)^*$,
and decomposition
$$L^2\Gamma({{\mathcal  H}^1_\cdot}^*)=(Ker \bar d^1\cap Ker d^*)
\oplus \overline{Image\  \bar d}
\oplus \overline {Image\  {{d^1}}^*},$$
we can see that it is self-adjoint.  By Von-Neumann's theorem
 \cite{Reed-Simon-vol2} Theorem X25, $(\tilde d+\tilde d^*)^2$ is also
self-adjoint (and in particular has dense domain). From Proposition
 \ref{pr:d2=0}
we see that for $(f,\theta, \phi)$ in its domain
$$(\tilde d+\tilde d^*)^2(f,\theta, \phi)=(\bar d^*d f ,
  ({\bar {d^1}}^* \bar d^1+\bar d\bar d^*)\theta,
 \bar d^1{\bar {d^1}}^*\phi)$$
and $$Dom(\tilde d+\tilde d^*)^2   =Dom(\bar d^*d)\oplus 
  Dom({\bar {d^1}}^*\bar d^1+\bar d\bar d^*)
    \oplus Dom (\bar d^1{\bar {d^1}}^*).$$
In particular if we set
 $$\Delta^1:={\bar {d^1}}^* \bar d^1+ \bar d \bar d^*$$
we obtain a nonnegative self-adjoint operator on
 $L\Gamma({{\mathcal  H}_\cdot^1}^*)$.
Since $$Ker(\tilde d+\tilde d^*)^2\ =\ Ker (\tilde d+\tilde d^*)\ = \ 
Ker\ \tilde d\cap Ker \ \tilde d^*$$ we see that
 $\phi\in L^2\Gamma({{\mathcal  H}^1_\cdot}^*)$ is harmonic,
 i.e. $\phi\in Ker \Delta^1$,
if and only if $\bar d^1\phi=0$ and $\bar d^*\phi=0$.

\bigskip

As remarked in \S\ref{section-exterior} we know $\bar d$ has closed
range.  Thus:
\begin{theorem}
The space $L^2\Gamma({{\mathcal  H}^1}^*)$ of ${\mathcal  H}$ 1-forms has the
 decomposition
$$L^2\Gamma({{\mathcal  H}^1_\cdot}^*)=Ker\Delta^1\oplus Image\  \bar d
\oplus \overline {Image\  { {d^1}}^*}.$$
In particular every cohomology class in $L^2H^1(C_{x_0}M)$ has a
unique representative in $Ker\Delta^1$.
\end{theorem}

\end{document}